\documentclass[a4paper,12pt]{article}
\usepackage[utf8]{inputenc}

\usepackage{amsmath}
\usepackage{amsthm}
\usepackage{amssymb}
\usepackage{mathtools}
\usepackage{authblk}

\def\R{\mathbb{R}}
\def\C{\mathbb C}
\def\Z{\mathbb Z}
\def\N{\mathbb N}
\def\P{\mathbb P}

\newtheorem*{remark*}{Remark}
\newtheorem{example}{Example}
\newtheorem{theorem}{Theorem} 
\newtheorem{lemma}{Lemma}

\title{THE LIMIT CYCLES IN A GENERALIZED RAYLEIGH-LI\'ENARD OSCILLATOR}
\author[a]{Lubomir Gavrilov}
\author[b]{Iliya D. Iliev}
 \affil[a]{ Institut de Math\'{e}matiques de Toulouse ; UMR 5219
 
  Universit\'{e} de Toulouse ; CNRS 
           
    UPS, F-31062 Toulouse Cedex 9,  France}
\affil[b]{ Institute of Mathematics, Bulgarian Academy of Sciences

 Bl. 8, 1113 Sofia, Bulgaria}

\begin{document}

\maketitle

\begin{abstract}
We compute the cyclicity of open period annuli of the following generalized Rayleigh-Li\'enard equation
$$\ddot{x}+ax+bx^3-(\lambda_1+\lambda_2 x^2+\lambda_3\dot{x}^2+\lambda_4  x^4+\lambda_5\dot{x}^4+\lambda_6 x^6)\dot{x}=0$$
and the equivalent planar system $X_\lambda$,
where the coefficients of the perturbation $\lambda_j$ are independent small parameters and $a, b$ are fixed nonzero constants.
Our main tool is the machinery of the so called higher-order Poincar\'e-Pontryagin-Melnikov functions (Melnikov functions $M_n$ for short), combined with the explicit computation of center conditions  and the corresponding Bautin ideal. 

We consider first arbitrary analytic arcs $\varepsilon \to \lambda(\varepsilon)$ and explicitly compute all possible Melnikov functions $M_n$ related to the deformation $X_{ \lambda(\varepsilon)} $. At a second step we obtain exact bounds for the number of the zeros of the Melnikov functions (complete elliptic integrals depending on parameter)  in an appropriate complex domain, using a modification of Petrov's method. 

To deal with the general case of six-parameter deformations $\lambda \to X_\lambda$, we compute first the related Bautin ideal. To do this we carefully study the  Melnikov functions up to order three, and then use Nakayama lemma from Algebraic geometry. The principalization of the Bautin ideal (achieved after a blow up) reduces finally the study of general deformations $X_{ \lambda } $ to the study of one-parameter deformations $X_{ \lambda(\varepsilon)} $.
\end{abstract}

\tableofcontents


\section{Statement of the results}
We study the limit cycles of the following equations
\begin{align}
\label{eq1}
\ddot{x}+ax+bx^3-(\lambda_1+\lambda_2 x^2+\lambda_3\dot{x}^2+\lambda_4  x^4+\lambda_5\dot{x}^4+\lambda_6 x^6)\dot{x}=0
\end{align}
assuming that $a$ and $b$ are fixed non-zero constants, and $\lambda_j$ are small real parameters.

When  $\lambda_j$ are all equal to zero, the system has a polynomial first integral $$H=\frac12 \dot{x}^2+\frac12ax^2+\frac14bx^4 .
$$

When $\lambda_3
=\lambda_5=0$, (\ref{eq1}) is Li\'enard equation, and if $\lambda_2=\lambda_4=\lambda_6=0$,(\ref{eq1}) is Rayleigh equation. 
For historical comments  on this generalized Rayleigh-Li\'enard equation
(\ref{eq1}) we refer to the recent paper \cite{ELT}.

Equivalently, taking 
$$H=\frac12 y^2+\frac12ax^2+\frac14bx^4, \;\lambda=(\lambda_1, \lambda_2, ..., \lambda_6)$$
one obtains the equivalent planar differential system 
(vector field $X_\lambda$)
\begin{align}
\label{eq2}
\tag{$2^\lambda$}
X_\lambda : 
\left\{ 
\begin{array}{ll}\dot{x}=H_y,\\[2mm] 
\dot{y}=-H_x+ (\lambda_1+\lambda_2 x^2+\lambda_3y^2+\lambda_4  x^4+\lambda_5y^4+\lambda_6 x^6)y,
\end{array}
\right.
\end{align}
which will be considered as a small deformation of the Hamiltonian system
\begin{align}
\label{eq20}
\tag{$2^0$}
X_0 :  \left\{ 
\begin{array}{ll}\dot{x}=H_y,\\[2mm] 
\dot{y}=-H_x,
\end{array}
\right.
\end{align}

Our work was motivated by \cite{WHC} and the recent paper \cite{ELT} which study  bifurcations of limit cycles  in {\rm(\ref{eq2})}, provided that all $\lambda_j$ depend linearly on a small parameter $\varepsilon$. In fact, \cite{ELT} deals only with  
limit cycles bifurcating from separatrix connections of (\ref{eq20}).

\begin{figure}
\begin{center}
\includegraphics[width=12cm]{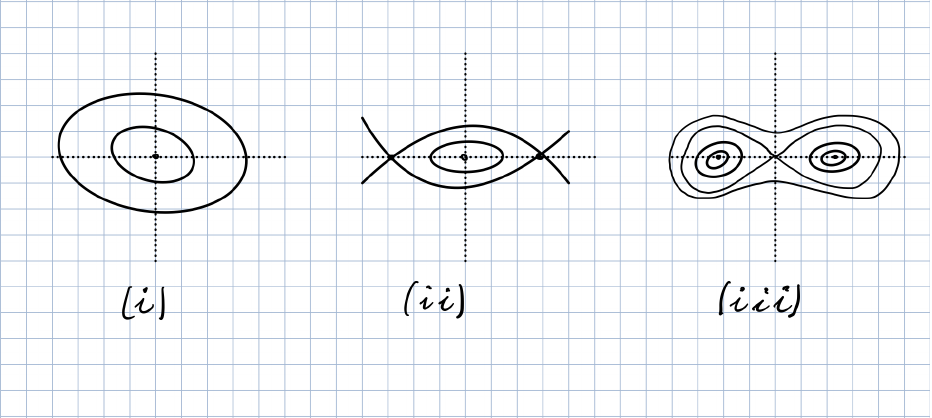}
\end{center}
\caption{ Phase portraits of ($2^0$) : the global center, the truncated pendulum and the eight loop}
\label{figure1a}
\end{figure}

Clearly, by scaling the variables, one can assume without loss of generality that $|a|=|b|=1$. Moreover, in the case of negative $a,b$ and arbitrary $\lambda_j$ the system (\ref{eq2})
has a single equilibrium point which is a saddle, hence it has no limit cycles. 
Then the following three cases of interest appear, see Fig.\ref{figure1a} :
$$\begin{array}{ll}
(i) & H=\frac12y^2+\frac12x^2+\frac14x^4, \;\mbox{\rm (global center)} \quad h>0;\\[2mm]
(ii) & H=\frac12y^2+\frac12x^2-\frac14x^4, \;\mbox{\rm (truncated pendulum)} \quad 0<h<\frac14;\\[2mm]
(iii) & H=\frac12y^2-\frac12x^2+\frac14x^4, \;\mbox{\rm (eight loop)} \quad -\frac14<h<0\cup h>0;
\end{array}$$
where the Hamiltonian level $h$ runs the intervals forming a set of nested ovals in the phase space (period annulus of the Hamiltonian system). 
Case (iii) is also known as the {\it Duffing oscillator}. 

In the next four sub-sections we describe the main results of the paper.

\subsection{Limit cycles} 
We first recall that the cyclicity of an annulus $\Pi$  is the maximal number of limit cycles which tend to the compact subsets $K\subset \Pi$ 
when $\lambda$ tends to zero \cite{Roussarie,gavr08}. The main result of the paper is the following
\begin{theorem}  
\label{main}
The cyclicity of each given open period annulus of the Hamiltonian system {\rm (\ref{eq20})} with respect to the six-parameter deformation  
{\rm(\ref{eq2})} is five, except for the exterior period annulus of the eight loop case. In this latter case, the cyclicity is bounded by six.
\end{theorem}
Thus, except for the exterior period annulus of the eight loop case, the above Theorem gives an  exact answer for the cyclicity of open period annuli and can not be improved.
The exterior period annulus of the eight loop case is the only one which does not contain a center in its closure. 

\emph{We conjecture, that the cyclicity of this exterior open period annulus is exactly $6$.}

\begin{remark*}
\begin{enumerate} \mbox{ }
\item 
In the case when a period annulus of {\rm(\ref{eq20})} contains in its closure an  equilibrium point of center type $C$,  denote  $\tilde{\Pi} = \Pi\cup C$. Then the cyclicity of 
$\tilde{\Pi}$ is also five.  

\item
 By a well known result of  Roussarie \cite[Theorem 25]{Roussarie} it is possible to estimate the cyclicity of homoclinic loops, which would imply also estimates for the cyclicity of bounded period annuli.
\item
The vector field (\ref{eq2}) has central symmetry with respect to the origin. Therefore in the eight loop case (iii) the distribution of limit cycles bifurcating from 
the bounded period annuli is the same. This implies that the total cyclicity of the union of the bounded period annuli in the eight loop case (iii) equals 10.
\item
On the other hand, in 
the  global center case case (i) and in the eight loop case (iii) we can say nothing about limit cycles bifurcating from infinity, nor in the truncated pendulum case (ii) about limit cycles bifurcating from the connection containing $\Pi$ in its interior. In the eight loop case (iii) Theorem \ref{main} says nothing  about the cyclicity of the closure of the unbounded period annulus.
\end{enumerate}
\end{remark*}

In the remaining three sub-sections we indicate the intermediate results (some of them of independent interest) needed for the proof of the Theorem\ref{main}.

\subsection{Melnikov functions related to analytic arcs}  To prove the above result we consider initially analytic arcs 
\label{section12}
$$\varepsilon \to \lambda(\varepsilon) , \; \lambda(0)=0$$
together with the corresponding one-parameter deformations $\varepsilon \to X_{\lambda(\varepsilon)}$, and then use the powerful machinery of the so called higher order Poincar\'e-Pontryagin-Melnikov functions, see \cite{JPF, gavr08, gavr20}. Namely, let $d(h,\varepsilon)$ be the displacement function associated to period annulus $\Pi$ of $X_0$ with respect to the deformation $X_{\lambda(\varepsilon)}$, where $h$ is the restriction of $H(x,y)$ on a cross section to the period orbits in $\Pi$. As $d(h,\varepsilon)$ is analytic both in $h$ and $\varepsilon$ then we have
\begin{align}
\label{eq3}
\stepcounter{equation}
d(h,\varepsilon)=M_1(h)\varepsilon+M_2(h)\varepsilon^2+M_3(h)\varepsilon^3+\ldots
\end{align}
and let $M_n$ be the first non-vanishing coefficient in this series. Then $M_n$ is the Poincar\'e-Pontryagin-Melnikov  (or bifurcation) function, and 
our primary goal will be its calculation. Indeed,  the zeros of the first non-vanishing coefficient $M_n(h)$, $n\in \N$ give the number, multiplicity and the positions of the periodic orbits in the period annulus of the original Hamiltonian system which persist under the perturbation $X_{\lambda(\varepsilon)}$
becoming limit cycles. 
Let us denote 
\begin{align}
I_0(h)=\int_{H=h}ydx, I_2(h)=\int_{H=h}x^2ydx
\end{align}
 where $h$ runs the respective interval. Below, we denote by $m$ a natural number. We shall prove the following:

\begin{theorem}
\label{th1}
Let $M_n(h)$ be the first nonzero coefficient in the expansion of the displacement function at
 $H=h$, namely $$d(h,\varepsilon)=M_n(h)\varepsilon^n + M_{n+1}(h)\varepsilon^{n+1}+\ldots, n\in \N .$$
Then $M_n(h)$ has for any of the three cases the form
\begin{align*}
M_n(h)&=(p_{2n}^1h^2+p_{1n}^1h+p_{0n}^1)I_2(h)+(q_{2n}^1h^2+q_{1n}^1h+q_{0n}^1)I_0(h),\quad &n=3m,\\
M_n(h)&=(p_{1n}^0h+p_{0n}^0)I_2(h)+(q_{2n}^0h^2+q_{1n}^0h+q_{0n}^0)I_0(h),\quad &n\neq 3m,
\end{align*}
with constants $p^i_{jn}$ and $q^i_{jn}$ explicitly determined by the coefficients in $(1)$.
\end{theorem}

Assuming that 
$$\lambda_j=\lambda_j(\varepsilon)=\sum_{k=1}^\infty \lambda_{jk}\varepsilon^k,\quad 1\leq j\leq 6$$

\vspace{1ex}
\noindent 
then the exact formulation for any of the cases is as follows: 

\begin{theorem}[{\small\sc global center}]
\label{th2}
 {\it Let $M_n(h)$ be the first nonzero coefficient in the expansion of the 
displacement function at $H=h$, $h>0$. Then for any} $n \in \N$ 
$$\begin{array}{l}M_n(h)= [(\frac43\lambda_{6n}-\frac{320}{231}\lambda_{5n})h+(\frac{32}{21}\lambda_{6n}-\frac{40}{231}\lambda_{5n}-\frac87\lambda_{4n}-\frac37\lambda_{3n}+\lambda_{2n})]I_2(h)\\[2mm]
   +[\frac{240}{77}\lambda_{5n}h^2+(-\frac{16}{21}\lambda_{6n}+\frac{20}{231}\lambda_{5n}+\frac47\lambda_{4n}
+\frac{12}{7}\lambda_{3n})h+\lambda_{1n}]I_0(h)\\[2mm]
+\frac{8}{1001}\lambda_{3m}^3[-(308h^2+181h+24)I_2(h)+(80h^2+12h)I_0(h)]\; \mbox{\it if}\; n=3m,\\[4mm]
M_n(h)= [(\frac43\lambda_{6n}-\frac{320}{231}\lambda_{5n})h+(\frac{32}{21}\lambda_{6n}-\frac{40}{231}\lambda_{5n}-\frac87\lambda_{4n}-\frac37\lambda_{3n}+\lambda_{2n})]I_2(h)\\[2mm]
   +[\frac{240}{77}\lambda_{5n}h^2+(-\frac{16}{21}\lambda_{6n}+\frac{20}{231}\lambda_{5n}+\frac47\lambda_{4n}
+\frac{12}{7}\lambda_{3n})h+\lambda_{1n}]I_0(h)\; \mbox{\it if}\; n\neq 3m.
\end{array}$$
\end{theorem}
\vspace{2ex}
\noindent
\begin{theorem}[{\small\sc truncated pendulum}]
\label{th3}
Let $M_n(h)$ be the first nonzero coefficient in the expansion of the displacement function at $H=h$, $h\in (0,\frac14)$. Then for any $n\in\N$  
$$\begin{array}{l}M_n(h)= [-(\frac43\lambda_{6n}+\frac{320}{231}\lambda_{5n})h+(\frac{32}{21}\lambda_{6n}+\frac{40}{231}\lambda_{5n}+\frac87\lambda_{4n}-\frac37\lambda_{3n}+\lambda_{2n})]I_2(h)\\[2mm]
   +[\frac{240}{77}\lambda_{5n}h^2-(\frac{16}{21}\lambda_{6n}+\frac{20}{231}\lambda_{5n}+\frac47\lambda_{4n}
-\frac{12}{7}\lambda_{3n})h+\lambda_{1n}]I_0(h)\\[2mm]
-\frac{8}{1001}\lambda_{3m}^3[(308h^2-181h+24)I_2(h)+(80h^2-12h)I_0(h)]\; \mbox{\it if}\; n=3m,\\[4mm]
M_n(h)=[-(\frac43\lambda_{6n}+\frac{320}{231}\lambda_{5n})h+(\frac{32}{21}\lambda_{6n}+\frac{40}{231}\lambda_{5n}+\frac87\lambda_{4n}-\frac37\lambda_{3n}+\lambda_{2n})]I_2(h)\\[2mm]
   +[\frac{240}{77}\lambda_{5n}h^2-(\frac{16}{21}\lambda_{6n}+\frac{20}{231}\lambda_{5n}+\frac47\lambda_{4n}
-\frac{12}{7}\lambda_{3n})h+\lambda_{1n}]I_0(h)\; \mbox{\it if}\; n\neq 3m.
\end{array}$$
\end{theorem}

\begin{theorem}[{\small\sc eight loop}]
\label{th4}
 Let $M_n(h)$ be the first nonzero coefficient in the expansion of the displacement 
function at $H=h$, $h\in (-\frac14,0)\cup(0,\infty)$. Then for any $n\in\N$  
$$\begin{array}{l}M_n(h)= [(\frac43\lambda_{6n}+\frac{320}{231}\lambda_{5n})h+(\frac{32}{21}\lambda_{6n}+\frac{40}{231}\lambda_{5n}+\frac87\lambda_{4n}+\frac37\lambda_{3n}+\lambda_{2n})]I_2(h)\\[2mm]
   +[\frac{240}{77}\lambda_{5n}h^2+(\frac{16}{21}\lambda_{6n}+\frac{20}{231}\lambda_{5n}+\frac47\lambda_{4n}
+\frac{12}{7}\lambda_{3n})h+\lambda_{1n}]I_0(h)\\[2mm]
-\frac{8}{1001}\lambda_{3m}^3[(308h^2+181h+24)I_2(h)+(80h^2+12h)I_0(h)]\; \mbox{\it if}\; n=3m,\\[4mm]
M_n(h)= [(\frac43\lambda_{6n}+\frac{320}{231}\lambda_{5n})h+(\frac{32}{21}\lambda_{6n}+\frac{40}{231}\lambda_{5n}+\frac87\lambda_{4n}+\frac37\lambda_{3n}+\lambda_{2n})]I_2(h)\\[2mm]
   +[\frac{240}{77}\lambda_{5n}h^2+(\frac{16}{21}\lambda_{6n}+\frac{20}{231}\lambda_{5n}+\frac47\lambda_{4n}
+\frac{12}{7}\lambda_{3n})h+\lambda_{1n}]I_0(h)\; \mbox{\it if}\; n\neq 3m.
\end{array}$$
\end{theorem}
\vspace{1ex}
\vspace{1ex}
\noindent
We would like to mention that quite similar results also hold if $x^6$ in the perturbation is replaced by $x^2y^2$
 (equivalently, by $x^2\dot{x}^2$). 

It is also worth noting that the change of variables transforming the original Hamiltonian to $(i), (ii)$ and $(iii)$ is
$t=\frac{1}{\sqrt{|a|}}t_1$, $x(t)=\sqrt{|\frac{a}{b}|}x_1(t_1)$, $y(t)=\frac{|a|}{\sqrt{|b|}}y_1(t_1)$. 
Therefore if someone is interested in the impact of the Hamiltonian parameters on $M_n(h)$, then it is easy to 
see that the present set of lambdas should be replaced by 

\begin{align}
\label{new}
\frac{1}{\sqrt{|a|b^6}}\times \{|b^3|\lambda_1, |a|b^2\lambda_2, 
a^2b^2\lambda_3, a^2|b|\lambda_4, a^4|b|\lambda_5, |a|^3\lambda_6\}.
\end{align}

\vspace{1ex}
Let us recall that for a general  planar Hamiltonian system perturbed by small perturbation which is analytic with respect 
to the small parameter $\varepsilon$, namely 
$$\begin{array}{ll}\dot{x}=H_y+f(x,y,\varepsilon), & f(x,y,\varepsilon)=\varepsilon f_1(x,y)+\varepsilon^2 f_2(x,y)
+\varepsilon^3 f_3(x,y)+\ldots,\\[2mm]
\dot{y}=-H_x+g(x,y,\varepsilon), & g(x,y,\varepsilon)=\varepsilon g_1(x,y)+\varepsilon^2 g_2(x,y)+\varepsilon^3 g_3(x,y)+\ldots,
\end{array}$$
the procedure of calculating $M_n(h)$ is the following, see e.g. \cite{IDI}, \cite{Roussarie} and the pioneering work by Fran\c{c}oise \cite{JPF}. Fix an open period annulus $\Pi= \{(x,y)\}: H(x,y)=h\in (\alpha,\beta)\subset \R \}$. 
Take the one-forms  $\omega_k=g_k(x,y)dx-f_k(x,y)dy$, $k\in \N$, $\omega(\varepsilon)=\varepsilon\omega_1+
\varepsilon^2\omega_2+\varepsilon^3\omega_3+....$. Then the system written in Pfaffian form is simply $dH=\omega(\varepsilon)$.
When parametrized by the Hamiltonian level $h$, the displacement function developed in series is   (\ref{eq3})
where $M_k(h)=\oint_{H=h}\Omega_k$  with  $\Omega_1=\omega_1$, and provided that $\Omega_j$ is relatively exact form if $j\leq k$, that is
$\Omega_j=r_jdH+dR_j$, $j\leq k$, then $M_{k+1}(h)=\oint_{H=h}\Omega_{k+1}$, with $\Omega_{k+1}=\omega_{k+1}
+\sum_{i+j=k+1}r_j\omega_i$.
In particular, $\Omega_2=\omega_2+r_1\omega_1$, $\Omega_3=\omega_3+r_1\omega_2+r_2\omega_1$,
 $\Omega_4=\omega_4+r_1\omega_3+r_2\omega_2+r_3\omega_1$, and so on. Clearly, $\Omega_j$ is 
relatively exact form if and only if $\int_{H=h}\Omega_j\equiv 0$ in $\Pi$. This procedure was initially discovered by Fran\c{c}oise  
for the Hamiltonian $H=\frac12x^2+\frac12y^2$ and the perturbation $\varepsilon\omega_1$. 

To get some idea about the proof of Theorem \ref{th1}, we initially calculate the first 6 coefficients in the expression of the displacement function. Then we proceed to finish the proof by induction. Our main technical tool is the relative decomposition of monomial one-forms with respect to given a Hamiltonian $H$ and its period annulus $\Pi$ (if not unique). For this purpose, we prepared as in \cite{ILY} a list of decompositions we will need, see the Appendix. In our case the decompositions do not depend on $\Pi$ and their right-hand sides have the form 
$[u(H)x^2+v(H)]ydx+rdH+dR$. Once derived, these formulas would be very useful to calculate $M_k(h)$ as well as $r_k$ we need. Moreover, to calculate $M_k(h)$ itself, we use the first part of the triad on the right in each formula. By integrating it, we obtain immediately 
$[u(h)I_2(h)+v(h)I_0(h)]$ as a result. When the first part is missing, the coefficient $r$ at $dH$ is used to calculate $r_k$, a function we use in the next step of the construction. The last part $dR$ of each formula can be used to check it and in addition to produce new formulas by multiplying a known decomposition with appropriate monomial. All formulas in the Appendix could be verified by direct calculations. 

\subsection{Zeros of elliptic integrals}  Motivated by Theorem \ref{th1} we consider the vector space of functions
$$
V= Span\{ h^iI_0(h), h^jI_2(h): \; i,j = 0,1, 2 \}
$$
where $h$ belongs to one of the real intervals 
$$\textstyle (i) \; h>0,\; (ii)\; 0<h< \frac 14,\;(iii) -\frac 14 < h< 0\; {\rm or}\; h>0,$$ 
according to the case into consideration. Every Melnikov function belongs to $V$. A careful inspection of the expressions for the Melnikov functions, see   
Theorems  \ref{th2},  \ref{th3}, \ref{th4}, shows that the opposite is also true: every function in $V$ can be realized as  a Melnikov function. This can be also seen from the Bautin ideal, in the next section. The set of Melnikov functions is in one-to-one correspondence to the points on the exceptional divisor of the blowup of the Bautin ideal at $\lambda= 0$ \cite{gavr20}.
In our case the exceptional divisor is computed to be $\mathbb P^5$.

We shall prove

\begin{theorem}
\label{th6}
Except in the case  $0 < h \; (iii)$, 
the vector space of functions $V$ is Chebyshev, that is to say each function $I\in V$ can have at most $5 = \dim V - 1$ zeros in its interval of definition, where zeros are counted with multiplicity. In the case $0 < h \; (iii)$, the number of the zeros of each function $I\in V$ is bounded by $6$.
\end{theorem}
We shall prove in fact a stronger result, but for the analytic continuations of the derivatives  $I'(h)$ in an appropriate complex domain containing the real interval of definition.  
To evaluate the number of the zeros we use the argument principle, inspired by Petrov
 \cite{petr89,petr97}. The above Theorem, however, does not follow from these papers. Its proof contains moreover a new technical argument (Lemma \ref{lemma4}).
This implies also a stronger result which includes  the closure of the real interval definition, as it will be explained later in the text.

The above considerations already prove a bound for the cyclicity of  one-parameter deformations $X_{\lambda(\varepsilon)  }$. 


\subsection{The Bautin ideal}  To complete the proof of Theorem \ref{main} we need  the Bautin ideal of the annulus, associated to the general deformation $X_\lambda$. Denote the corresponding ideal, localized at $\lambda=0$, by $\mathcal B$.
Clearly $\mathcal B$ is an ideal of the ring of convergent power series $\R \{\lambda\}$ at $\lambda=0$. It is however polynomially generated and its zero locus (the center set) is just the origin $\{\lambda=0 \}$ (by Theorem \ref{th1}).
There are very few polynomial systems with explicitly known Bautin ideal. Therefore the next theorem is of independent interest.

\begin{theorem}
\label{th5ab}
Let $a, b$ be fixed non-zero constants, which are not both negative. The Bautin ideal $\mathcal B$ of the perturbed equations (\ref{eq2}) is given by
$$  \mathcal B = ( \lambda_{1},  \lambda_{2}+ 3a\lambda_{3}, \lambda_3^3, \lambda_{4} + 3b\lambda_{3}, \lambda_{5}, \lambda_{6})$$      
where, respectively, 
$$
\begin{array}{lll}
(i)      &\;\mbox{\rm (global center)} \;\; a>0, b>0, & \quad h>0; \\[2mm] 
(ii)     & \;\mbox{\rm (truncated pendulum)} \;\; a>0>b, & \quad 0<h<-\frac{a^2}{4b};\\[2mm]
(iii)    & \;\mbox{\rm (eight loop)} \;\; a<0<b,  & \quad -\frac{a^2}{4b}<h<0\cup h>0.
\end{array}
$$
\end{theorem}

After an obvious rescaling of $a, b$ we can suppose as in the beginning of the section, that $|a|=|b|=1$ and we are in one of the cases (i), (ii) or (iii). Then using (\ref{new})
we can equivalently restate
Theorem \ref{th5ab} as follows
\begin{theorem}
\label{th5}
The Bautin ideal $\mathcal B$ of the perturbed equation (\ref{eq2}) is given respectively by
$$\begin{array}{ll}
(i) & \mathcal B = ( \lambda_{1},    \lambda_{2}+ 3\lambda_{3} ,\lambda_3^3, \lambda_{4} + 3 \lambda_{3},    \lambda_{5},\lambda_{6})         , \\
&\;\mbox{\rm (global center)} \quad h>0;\\[2mm] 
(ii) & \mathcal B = ( \lambda_{1},    \lambda_{2}+ 3\lambda_{3} ,\lambda_3^3, \lambda_{4} - 3 \lambda_{3},    \lambda_{5},\lambda_{6})         , \\
& \;\mbox{\rm (truncated pendulum)} \quad 0<h<\frac14;\\[2mm]
(iii) & \mathcal B = ( \lambda_{1},    \lambda_{2} - 3\lambda_{3} ,\lambda_3^3, \lambda_{4} + 3 \lambda_{3},    \lambda_{5},\lambda_{6})         , \\
&   \;\mbox{\rm (eight loop)} \quad -\frac14<h<0\cup h>0;
\end{array}$$
\end{theorem}

The proof of Theorem \ref{th5} itself combines  the expansions found in Theorem \ref{th2},  \ref{th3} and   \ref{th4} with a version of Nakayama lemma in Algebraic Geometry, 
and might be of independent interest too.

\section{Calculation of the Melnikov functions in the global center case}

We begin by calculating the first 6 coefficients in the expression of the displacement function corresponding to the period annulus in the case of global center $(i)$ above.Thus, having in mind that $M_1(h)=\int_{H=h}\Omega_1$ where
 $$ \Omega_1=\omega_1=(\lambda_{11}+\lambda_{21} x^2+\lambda_{31}y^2+\lambda_{41}  x^4
+\lambda_{51}y^4+\lambda_{61} x^6)ydx$$
and using decompositions $(i1)-(i4)$  in the Appendix, we easily obtain
$$\begin{array}{rl}
M_1(h)= & \lambda_{11}I_0+\lambda_{21}I_2+\lambda_{31}(\frac{12}{7}hI_0-\frac37I_2)+\lambda_{41}(\frac47hI_0-\frac87I_2)\\[2mm]
           &+\lambda_{51}(\frac{240}{77}h^2I_0-\frac{320}{231}hI_2-\frac{40}{231}I_2+\frac{20}{231}hI_0)\\
           &+\lambda_{61}(\frac43hI_2
+\frac{32}{21}I_2-\frac{16}{21}hI_0)\\[2mm]
 =& [(\frac43\lambda_{61}-\frac{320}{231}\lambda_{51})h+(\frac{32}{21}\lambda_{61}-\frac{40}{231}\lambda_{51}-\frac87\lambda_{41}-\frac37\lambda_{31}+\lambda_{21})]I_2(h)\\[2mm]
   &+[\frac{240}{77}\lambda_{51}h^2+(-\frac{16}{21}\lambda_{61}+\frac{20}{231}\lambda_{51}+\frac47\lambda_{41}
+\frac{12}{7}\lambda_{31})h+\lambda_{11}]I_0(h).
\end{array}$$
Next, we see that $M_1(h)\equiv 0$ is equivalent to 
$\lambda_{11}=\lambda_{51}=\lambda_{61}=\lambda_{21}+3\lambda_{31}=\lambda_{41}+3\lambda_{31}=0$. Therefore $\omega_1$
reduces to $\lambda_{31}(-3x^2+y^2-3x^4)ydx=\lambda_{31}[-3xydH+d(xy^3)]$. In particular, $r_1=-3xy\lambda_{31}$. To handle $M_k(h)$, $k\geq 2$, we will use

\vspace{2ex}
\begin{lemma}
\label{lemma1}
({\small\sc global center})\\
  (a) Assume that $\omega_j=\lambda_{3j}[-3xydH+d(xy^3)]$ and also
 $r_k=-3xy\lambda_{3k}$. Then 
 $$r_k\omega_j=\lambda_{3j}\lambda_{3k}(8x^2y^2-4Hx^2-x^2-y^2)dH+ \mbox{\it exact form} .$$

\vspace{2ex}
\noindent
(b) {\it Let $\omega_j=\lambda_{3j}[-3xydH+d(xy^3)]$ and  $r_k=\Lambda_k(8x^2y^2-4Hx^2-x^2-y^2)-3xy\lambda_{3k}$ 
with $\Lambda_k=\sum_{l+m=k}\lambda_{3l}\lambda_{3m}$, $k\geq 2$. Then} 
$$\begin{array}{rl}
r_k\omega_j= & \lambda_{3j}\Lambda_k(-\frac{32}{5}x^2y^5+8Hx^2y^3+2x^2y^3-\frac25y^5)dx\\[2mm]
 &+\lambda_{3j}\Lambda_k(12Hx^3y-20x^3y^3+3x^3y+3xy^3)dH\\[2mm]
&+\lambda_{3j}\lambda_{3k}(8x^2y^2-4Hx^2-x^2-y^2)dH + \mbox{\it exact form}.
\end{array}$$
\end{lemma}

\vspace{1ex}
\noindent
{\bf Proof.} (a) We will perform the calculations modulo exact forms because we do not need their explicit expressions in our analysis. 
Thus we have 
$$\begin{array}{l}-3xy[-3xydH+d(xy^3)]=9x^2y^2dH-3xy(y^3dx+3xy^2dy)\\[2mm]
=9x^2y^2dH-3xy^4dx-\frac94x^2dy^4=9x^2y^2dH+\frac32xy^4dx+ \mbox{\it exact form}.\end{array}$$
Then using formula (i7) from the Appendix  yields the result in (a).

(b). The last part of the formula is the same as in (a). To obtain the first two expressions, we use the same calculation.
$$\begin{array}{l}
(8x^2y^2-4Hx^2-x^2-y^2)[-3xydH+d(xy^3)]\\[2mm]
=(-24x^3y^3+12Hx^3y+3x^3y+3xy^3)dH\\
+(8x^2y^5-4Hx^2y^3-x^2y^3-y^5)dx\\[2mm]
+(\frac{24}{5}x^3dy^5-4Hx^3dy^3-x^3dy^3-\frac35xdy^5).
\end{array}$$
The last expression in the brackets equals 
$$\textstyle 4x^3y^3dH+(-\frac{72}{5}x^2y^5+12Hx^2y^3+3x^2y^3+\frac35y^5)dx+ \mbox{\it exact form}.$$ 
which proves (b).$\Box$

\vspace{2ex}
\noindent
Now we are prepared to calculate $M_2(h)$ and $M_3(h)$ by integrating $\Omega_2=\omega_2+r_1\omega_1$ and 
 $\Omega_3=\omega_3+r_1\omega_2+r_2\omega_1$ respectively. Making use of Lemma 1(a)  with $j=k=1$, we see that
$r_1\omega_1$ is relatively exact form. Hence  $M_2(h)=\int_{H=h}\omega_2$ is just $M_1(h)$ with $\lambda_{j1}$ 
replaced by $\lambda_{j2}$:
$$\begin{array}{rl}
M_2(h)=& [(\frac43\lambda_{62}-\frac{320}{231}\lambda_{52})h+(\frac{32}{21}\lambda_{62}-\frac{40}{231}\lambda_{52}
-\frac87\lambda_{42}-\frac37\lambda_{32}+\lambda_{22})]I_2(h)\\[2mm]
   &+[\frac{240}{77}\lambda_{52}h^2+(-\frac{16}{21}\lambda_{62}+\frac{20}{231}\lambda_{52}+\frac47\lambda_{42}
+\frac{12}{7}\lambda_{32})h+\lambda_{12}]I_0(h).
\end{array}$$
Next, vanishing of $M_2(h)$ reduces $\omega_2$ to $\omega_2=\lambda_{32}[-3xydH+d(xy^3)]$ and therefore one can apply 
Lemma 1(a) again with $j=k=1$ to calculate first 
$$r_2=\lambda_{31}^2(8x^2y^2-4Hx^2-x^2-y^2)-3\lambda_{32}xy$$
and then with $j=2$, $k=1$ to see that $r_1\omega_2$ is relatively exact. Therefore 
$\int_{H=h}\Omega_3=\int_{H=h}(\omega_3+r_2\omega_1)$. The first integral yields an expression like $M_1(h)$ with all 
$\lambda_{j1}$ replaced by $\lambda_{j3}$. By Lemma 1(b) with $j=1, k=2$, the second integrand reduces to 
$\lambda^3_{31}(-\frac{32}{5}x^2y^5+8Hx^2y^3+2x^2y^3-\frac25y^5)dx$. 
Then we use formulas $(i3), (i5), (i6)$ from the Appendix which reduce the calculation to 
integrating the one-form (times $\lambda_{31}^3$)
$$\begin{array}{l}-\frac{32}{5}(\frac{80}{39}H^2x^2+\frac{3620}{3003}Hx^2-\frac{1600}{3003}H^2-\frac{80}{1001}H+\frac{160}{1001}x^2)ydx
\\[2mm]
+(8H+2)(\frac43Hx^2+\frac{8}{21}x^2-\frac{4}{21}H)ydx\\
-\frac25[\frac{240}{77}H^2-\frac{320}{231}Hx^2-\frac{40}{231}x^2+\frac{20}{231}H]ydx
\end{array}$$
which, after simplifying, leads to the final formula of $M_3(h)$ in terms of polynomial envelope of $I_0(h)$ and $I_2(h)$:
$$\begin{array}{rl}
M_3(h)=& [(\frac43\lambda_{63}-\frac{320}{231}\lambda_{53})h+(\frac{32}{21}\lambda_{63}-\frac{40}{231}\lambda_{53}
-\frac87\lambda_{43}-\frac37\lambda_{33}+\lambda_{23})]I_2(h)\\[2mm]
   &+[\frac{240}{77}\lambda_{53}h^2+(-\frac{16}{21}\lambda_{63}+\frac{20}{231}\lambda_{53}+\frac47\lambda_{43}
+\frac{12}{7}\lambda_{33})h+\lambda_{13}]I_0(h)\\[2mm]
&+\frac{8}{1001}\lambda_{31}^3[-(308h^2+181h+24)I_2(h)+(80h^2+12h)I_0(h)].
\end{array}$$
Obviously, $M_3(h)=0$ forces $\lambda_{31}$ be zero, together with $\omega_1$ and $r_1$. To summarize the results obtained so far, we formulate


\vspace{2ex}
\noindent
{\bf Corollary 1.} {\it Let $M_k(h), k\leq 3$ vanish. Then} 

\vspace{1ex}
\noindent
(i) $\lambda_{31}$, $\omega_1$, $r_1$ {\it are zero}, $\Lambda_2=\Lambda_3=0$, $\Lambda_4=\lambda_{32}^2$;

\vspace{1ex}
\noindent
(ii)  $\omega_k=\lambda_{3k}[-3xydH+d(xy^3)]$, $r_k=-3\lambda_{3k}xy$, $k=2,3$;

\vspace{1ex}
\noindent
{\bf Remark.} The simplest {\it essential perturbation} (see \cite[Definition 5]{gavr20}) with $M_1(h)=M_2(h)=0$ and $M_3(h)$ as above is
$$ Q(x,y,\varepsilon)=\varepsilon \lambda_{31}(-3x^2+y^2-3x^4)y+\varepsilon^3(\lambda_{13}+\lambda_{23}x^2
+\lambda_{43}x^4+\lambda_{53}y^4+\lambda_{63}x^6)y.$$

\vspace{1ex}
\noindent
Let us consider now the next triple set $M_k(h)=\int_{H=h}\Omega_k$, $k=4,5,6$ where (according to Corollary 1) 
$\Omega_4=\omega_4+r_2\omega_2$,  $\Omega_5=\omega_5+r_2\omega_3+r_3\omega_2$, 
$\Omega_6=\omega_6+r_2\omega_4+r_3\omega_3+r_4\omega_2$.
Applying Lemma 1(a) with $j=k=2$ we obtain as above that $M_4(h)=\int_{H=h}\omega_4$ takes the required form
 and moreover, when $M_4(h)$ vanishes,
$$\omega_4=\lambda_{34}[-3xydH+d(xy^3)],  r_4=\lambda_{32}^2(8x^2y^2-4Hx^2-x^2-y^2)-3\lambda_{34}xy.$$
Similarly, using twice  Lemma 1(a) with $j=2$, $k=3$ and $j=3$, $k=2$ one obtains that $M_5(h)=\int_{H=h}\omega_5$
is as in Theorem \ref{th2} and if $M_5(h)$ vanishes, then
$$\omega_5=\lambda_{35}[-3xydH+d(xy^3)],  r_5=2\lambda_{32}\lambda_{33}(8x^2y^2-4Hx^2-x^2-y^2)-3\lambda_{35}xy.$$
To perform the third step, we first use Lemma 1(a) twice with $j=4, k=2$ and $j=k=3$ to verify that $r_2\omega_4+r_3\omega_3$
is relatively exact form. Therefore $M_6(h)=\int_{H=h}(\omega_6+r_4\omega_2)$. To handle the second term, we use 
Lemma 1(b)
with $j=2, k=4$ to reduce it to $\lambda^3_{32}(-\frac{32}{5}x^2y^5+8Hx^2y^3+2x^2y^3-\frac25y^5)dx$ 
(modulo relatively exact form). In this way we obtain a similar formula of $M_6(h)$ like $M_3(h)$ above:
$$\begin{array}{rl}
M_6(h)=& [(\frac43\lambda_{66}-\frac{320}{231}\lambda_{56})h+(\frac{32}{21}\lambda_{66}-\frac{40}{231}\lambda_{56}
-\frac87\lambda_{46}-\frac37\lambda_{36}+\lambda_{26})]I_2(h)\\[2mm]
   &+[\frac{240}{77}\lambda_{56}h^2+(-\frac{16}{21}\lambda_{66}+\frac{20}{231}\lambda_{56}+\frac47\lambda_{46}
+\frac{12}{7}\lambda_{36})h+\lambda_{16}]I_0(h)\\[2mm]
&+\frac{8}{1001}\lambda_{32}^3[-(308h^2+181h+24)I_2(h)+(80h^2+12h)I_0(h)].
\end{array}$$
In particular if $M_6(h)$ vanishes, then $\lambda_{32}$ becomes zero, together with $\omega_2$ and $r_2$. Therefore one can 
formulate 

\vspace{2ex}
\noindent
{\bf Corollary 2.} {\it Let $M_k(h), k\leq 6$ vanish. Then} 

\vspace{1ex}
\noindent
(i) $\lambda_{3k}$, $\omega_k$, $r_k$ {\it are zero for} $k=1,2$, $\Lambda_k=0$, $2\leq k\leq 5$; $\Lambda_6=\lambda_{33}^2$;

\vspace{1ex}
\noindent
(ii)  $\omega_k=\lambda_{3k}[-3xydH+d(xy^3)]$, $3\leq k\leq 6$;

\vspace{1ex}
\noindent
(iii) $r_k=-3\lambda_{3k}xy$, $3\leq k\leq 5$; $r_6=\lambda_{33}^2(8x^2y^2-4Hx^2-x^2-y^2)-3\lambda_{36}xy$;

\vspace{2ex}
\noindent
{\bf Proof of Theorem \ref{th2}.}
To use induction in proving Theorem \ref{th2}, we need the following {\it inductive hypothesis} suggested by the corollaries above:

\vspace{1ex}
\noindent
 (a) {\it The functions $M_n(h)$ have for any $n\leq N=3K$ the expression as in} Theorem \ref{th2}. 

\vspace{1ex}
\noindent
(b) {\it If all $M_n(h)$, $n\leq N$ do vanish, then:}

\vspace{1ex}
\noindent
 (i) $\lambda_{3k}$, $\omega_k$, $r_k$ {\it vanish}, $1\leq k \leq K$; $\Lambda_k=0$, $2\leq k\leq 2K+1$,  
     $\Lambda_{2K+2}=\lambda_{3,K+1}^2$;

\vspace{1ex}
\noindent
(ii)  $\omega_k=\lambda_{3k}[-3xydH+d(xy^3)]$, $K+1\leq k\leq 3K$;

\vspace{1ex}
\noindent
(iii) $r_k=-3\lambda_{3k}xy$, $K+1\leq k\leq 2K+1$;

 $r_k=\Lambda_k(8x^2y^2-4Hx^2-x^2-y^2)-3\lambda_{3k}xy,\,\; 2K+2\leq k\leq 3K$ {\it if} $K\geq 2$.

\vspace{2ex}
\noindent 
From Corollaries 1 and 2 we are sure that the hypothesis holds if $K=1,2$. Assuming it holds with  $K$ we have 
to prove it holds with $K+1$. From part (b)(i) we conclude that the functions $M_{N+1}, M_{N+2}, M_{N+3}$ 
will come from integrating respectively the one forms (taken with $k=K$)
$$\begin{array}{l}
\Omega_{3k+1}=\omega_{3k+1}   +(r_{k+1}\omega_{2k}+\ldots+r_{2k}\omega_{k+1}),\\[2mm]
\Omega_{3k+2}=\omega_{3k+2}+(r_{k+1}\omega_{2k+1}+\ldots+r_{2k+1}\omega_{k+1})\; \mbox{\rm if}\; M_{N+1}=0,  \\[2mm]
\Omega_{3k+3}=\omega_{3k+3}+(r_{k+1}\omega_{2k+2}+\ldots+r_{2k+1}\omega_{k+2})+r_{2k+2}\omega_{k+1}\; \mbox{\rm if}\; M_{N+2}=0.
\end{array}$$
Moreover, by part (b)(ii)(iii) of the hypothesis and Lemma 1(a), the expressions in the brackets are relatively exact forms.
Therefore $M_{N+1}(h)=\int_{H=h}\omega_{3K+1}$ coincides with the expression in Theorem 3 and if $M_{N+1}(h)$ vanishes
then (b)(ii) above holds with $k=N+1$ too. Similarly, $M_{N+2}(h)=\int_{H=h}\omega_{3K+2}$ coincides with the expression in Theorem 3 
and if $M_{N+2}(h)=0$ then (b)(ii) above holds also with $k=N+2.$ Finally, 
  $$M_{N+3}(h)=\int_{H=h}(\omega_{3K+3}+r_{2K+2}\omega_{K+1}).$$The first term in the integrand yields just the first part in the formula 
in Theorem 3. The second term in the integrand, thanks to
statements (b)(ii) of the hypothesis applied with $k=K+1$, (b)(iii) applied with $k=2K+2$, (b)(i) concerning $\Lambda_{2K+2}$
and finally, Lemma 1(b) applied with $j=K+1, k=2K+2$ all together, becomes
\begin{align*}
\textstyle r_{2K+2}\omega_{K+1} & = \lambda_{3,K+1}^3(-\frac{32}{5}x^2y^5+8Hx^2y^3+2x^2y^3-\frac25y^5)dx\\
&+ \mbox{\it relatively exact form}.
\end{align*}
All this proves part (a) of the hypothesis  for $n\leq N+3=3(K+1)$. It still remains to establish that  part (b) of the hypothesis also 
holds with $K\to K+1$. Assume that $M_{3K+3}(h)$ also vanishes. This immediately  implies that  $\lambda_{3, K+1}$ together with
$\omega_{K+1}, r_{K+1}$, $\Lambda_{2K+2}, \Lambda_{2K+3}$ vanish too, while $\Lambda_{2K+4}=\lambda_{3,K+2}^2$
and $\omega_{3K+3}=\lambda_{3K+3}[-3xydH+d(xy^3)]$.
Therefore (b)(i) holds with $K$ replaced by $K+1$. Next, as we already mentioned above, (b)(ii) holds with $K$ replaced by $K+1$.
It only remains to check (b)(iii). As we already proved that $r_{K+1}=0$ and $\Lambda_{2K+2}=\Lambda_{2K+3}=0$, one obtains that 
the first formula in (b)(iii) holds with $K$ replaced by $K+1$. As $r_{K+1}=\omega_{K+1}=0$, the expressions of $\Omega_{3K+1}, 
\Omega_{3K+2}, \Omega_{3K+3}$ above (taken with $k=K$) reduce to 
$$\begin{array}{l}
\Omega_{3k+1}=\omega_{3k+1}+(r_{k+2}\omega_{2k-1}+\ldots+r_{2k-1}\omega_{k+2}),\\[2mm]
\Omega_{3k+2}=\omega_{3k+2}+(r_{k+2}\omega_{2k}+\ldots+r_{2k}\omega_{k+2}),\\[2mm]
\Omega_{3k+3}=\omega_{3k+3}+(r_{k+2}\omega_{2k+1}+\ldots+r_{2k+1}\omega_{k+2}).
\end{array}$$
Therefore, using the facts already proved above and applying Lemma 1(a) to each term, we obtain (with $k=K$)
 $$\begin{array}{l} r_{3k+1}=\Lambda_{3k+1}(8x^2y^2-4Hx^2-x^2-y^2)-3\lambda_{3, 3k+1}xy,\\[2mm]
                           r_{3k+2}=\Lambda_{3k+2}(8x^2y^2-4Hx^2-x^2-y^2)-3\lambda_{3, 3k+2}xy,\\[2mm]
                           r_{3k+3}=\Lambda_{3k+3}(8x^2y^2-4Hx^2-x^2-y^2)-3\lambda_{3, 3k+3}xy,\\[2mm]
\end{array}$$
which means that the second formula in (b)(iii) holds with $K$ replaced by $K+1$. 
The induction step procedure and together the proof of Theorem \ref{th2} are completed. $\Box$

\vspace{2ex}
\section{Calculation of the Melnikov functions in the truncated pendulum and the eight loop cases}

Just as in Section 2, we can use formulas with $(ii)$ and $(iii)$ from the Appendix in order to obtain the expressions 
of $M_1(h)$ which coincide with these with $n=1$ in Theorem \ref{th3} and Theorem \ref{th4}, respectively. In particular, if $M_1(h)$ vanish,
the one form $\omega_1$ reduces to 
$$
\begin{array}{l}
\omega_1=\lambda_{31}(-3x^2+y^2+3x^4)ydx=\lambda_{31}[-3xydH+d(xy^3)] \\[2mm] 
\;\; \mbox{\rm (truncated pendulum),} \\[2mm] 
\omega_1=\lambda_{31}(3x^2+y^2-3x^4)ydx=\lambda_{31}[-3xydH+d(xy^3)] \\
\;\; \mbox{\rm (eight loop)}.\end{array}
$$ 
In particular, $r_1=-3xy\lambda_{31}$ in both cases. As above, to handle $M_k(h)$, $k\geq 2$, we  use


\vspace{2ex}
\begin{lemma}
\label{lemma2}
({\small\sc truncated pendulum}) \\
(a) Assume that $\omega_j=\lambda_{3j}[-3xydH+d(xy^3)]$ and 
 $r_k=-3xy\lambda_{3k}$. Then 
 $$r_k \omega_j=\lambda_{3j} \lambda_{3k}(8x^2y^2-4Hx^2+x^2+y^2)dH + \mbox{\it exact form} .$$

\vspace{2ex}
\noindent
(a)  Let $\omega_j=\lambda_{3j}[-3xydH+d(xy^3)]$ and  $r_k=\Lambda_k(8x^2y^2-4Hx^2+x^2+y^2)-3xy\lambda_{3k}$ 
with $\Lambda_k=\sum_{l+m=k}\lambda_{3l}\lambda_{3m}$, $k\geq 2$. Then
$$
\begin{array}{rl}
r_k\omega_j= & \lambda_{3j}\Lambda_k(-\frac{32}{5}x^2y^5+8Hx^2y^3-2x^2y^3+\frac25y^5)dx\\[2mm]
 &+\lambda_{3j}\Lambda_k(12Hx^3y-20x^3y^3-3x^3y-3xy^3)dH \\[2mm]
 &+\lambda_{3j}\lambda_{3k}(8x^2y^2-4Hx^2+x^2+y^2)dH + \mbox{\it exact form}.
\end{array}
$$
\end{lemma}

\vspace{1ex}
\noindent
{\bf Inductive hypothesis 2} ({\small\sc truncated pendulum}):

\vspace{1ex}
\noindent
 (a) {\it The functions $M_n(h)$ have for any $n\leq N=3K$ the expression as in} Theorem \ref{th3}. 

\vspace{1ex}
\noindent
(b) {\it If all $M_n(h)$, $n\leq N$ do vanish, then:}

\vspace{1ex}
\noindent
 (i) $\lambda_{3k}$, $\omega_k$, $r_k$ {\it vanish}, $1\leq k \leq K$; $\Lambda_k=0$, $2\leq k\leq 2K+1$,  
     $\Lambda_{2K+2}=\lambda_{3,K+1}^2$;

\vspace{1ex}
\noindent
(ii)  $\omega_k=\lambda_{3k}[-3xydH+d(xy^3)]$, $K+1\leq k\leq 3K$;

\vspace{1ex}
\noindent
(iii) $r_k=-3\lambda_{3k}xy$, $K+1\leq k\leq 2K+1$;

 $r_k=\Lambda_k(8x^2y^2-4Hx^2+x^2+y^2)-3\lambda_{3k}xy,\,\; 2K+2\leq k\leq 3K$ {\it if} $K\geq 2$.

\vspace{2ex}
\begin{lemma}
\label{lemma3}
({\small\sc eight loop})  (a) {\it Assume that $\omega_j=\lambda_{3j}[-3xydH+d(xy^3)]$ and 
 $r_k=-3xy\lambda_{3k}$. Then $r_k\omega_j=\lambda_{3j}\lambda_{3k}(8x^2y^2-4Hx^2-x^2+y^2)dH$ + exact form}.

\vspace{2ex}
\noindent
(b) {\it Let $\omega_j=\lambda_{3j}[-3xydH+d(xy^3)]$ and  $r_k=\Lambda_k(8x^2y^2-4Hx^2-x^2+y^2)-3xy\lambda_{3k}$ 
with $\Lambda_k=\sum_{l+m=k}\lambda_{3l}\lambda_{3m}$, $k\geq 2$. Then} 
$$\begin{array}{rl}
r_k\omega_j= & \lambda_{3j}\Lambda_k(-\frac{32}{5}x^2y^5+8Hx^2y^3+2x^2y^3+\frac25y^5)dx\\[2mm]
 &+\lambda_{3j}\Lambda_k(12Hx^3y-20x^3y^3+3x^3y-3xy^3)dH\\[2mm]
&+\lambda_{3j}\lambda_{3k}(8x^2y^2-4Hx^2-x^2+y^2)dH + \mbox{\it exact form}.
\end{array}$$
\end{lemma}

\vspace{1ex}
\noindent
{\bf Inductive hypothesis 3} ({\small\sc eight loop}):

\vspace{1ex}
\noindent
 (a) {\it The functions $M_n(h)$ have for any $n\leq N=3K$ the expression as in} Theorem \ref{th4}. 

\vspace{1ex}
\noindent
(b) {\it If all $M_n(h)$, $n\leq N$ do vanish, then:}

\vspace{1ex}
\noindent
 (i) $\lambda_{3k}$, $\omega_k$, $r_k$ {\it vanish}, $1\leq k \leq K$; $\Lambda_k=0$, $2\leq k\leq 2K+1$,  
     $\Lambda_{2K+2}=\lambda_{3,K+1}^2$;

\vspace{1ex}
\noindent
(ii)  $\omega_k=\lambda_{3k}[-3xydH+d(xy^3)]$, $K+1\leq k\leq 3K$;

\vspace{1ex}
\noindent
(iii) $r_k=-3\lambda_{3k}xy$, $K+1\leq k\leq 2K+1$;

 $r_k=\Lambda_k(8x^2y^2-4Hx^2-x^2+y^2)-3\lambda_{3k}xy,\,\; 2K+2\leq k\leq 3K$ {\it if} $K\geq 2$.

Then following exactly the same way as in Section 2, including induction procedure, we verify the statements in Theorem \ref{th3} 
and Theorem \ref{th4}.  For this reason, we will omit repeating the same details. This finishes the proof of Theorem \ref{th1}, too.


\section{Zeros of elliptic integrals}
Here we prove Theorem \ref{th6}.
Full details of the proof will be given only
 for the exterior period annulus of the so called eight-loop case (iii)
$$
 H=\frac12y^2-\frac12x^2+\frac14x^4 .
 $$
 Indeed, this turns out to be   the most difficult part in the proof of Theorem \ref{th6}. 
 
Denote $$I_0(h)=\int_{\gamma(h)}ydx, I_2(h)=\int_{\gamma(h)}x^2ydx$$
 where $h>0$, and the real oval $\gamma(h) =\{H=h\}, h>0$ of the complex elliptic curve 
$$E_h = \{(x,y)\in \C^2 : H(x,y)=h\}$$
represents a homology class of $H_1(E_h,\Z)$ which by abuse of notation will be denoted by $\gamma(h)$ too.
We are interested in the real zeros, counted with multiplicity, of the complete elliptic integral of the form
$$
I(h)= p(h) I_2(h) + q(h) I_0(h)
$$
where $p,q$ are real polynomials of degree at most two, on the real interval $\{h : h>0\}$.

Instead of this we shall count the zeros of its derivative $I'(h)$ in an appropriate complex domain, containing $\{h : h>0\}$.
Denote
$$
J(h)=I'(h), J_2(h)=I_2'(h), J_0(h)=I_0'(h) .
$$
Using the Picard-Fuchs equation satisfied by $I_0, I_2$, see \cite[Lemma 3.1]{ ilpe99}
\begin{align*}
4h I_0' + I_2' & = 3 I_0 \\
4hI_0' + (12h+4) I_2' & = 15I_2
\end{align*}
we may express $J(h)=I'(h)$ in terms of $J_0, J_2$. The result is that
$$
J(h)= I'(h)=  \tilde p(h) J_2(h) + \tilde q(h) J_0(h)
$$
where $\tilde p, \tilde q$ might be  any real polynomials of degree at most two again. The continuous family of cycles $\gamma(h) \in H_1(E_h, \Z), h>0$,  is extended to a continuous family on the complex domain
$\mathcal D = \{ \C \setminus (-\infty, 0] \}$
\begin{figure}
\begin{center}
\includegraphics[width=10cm]{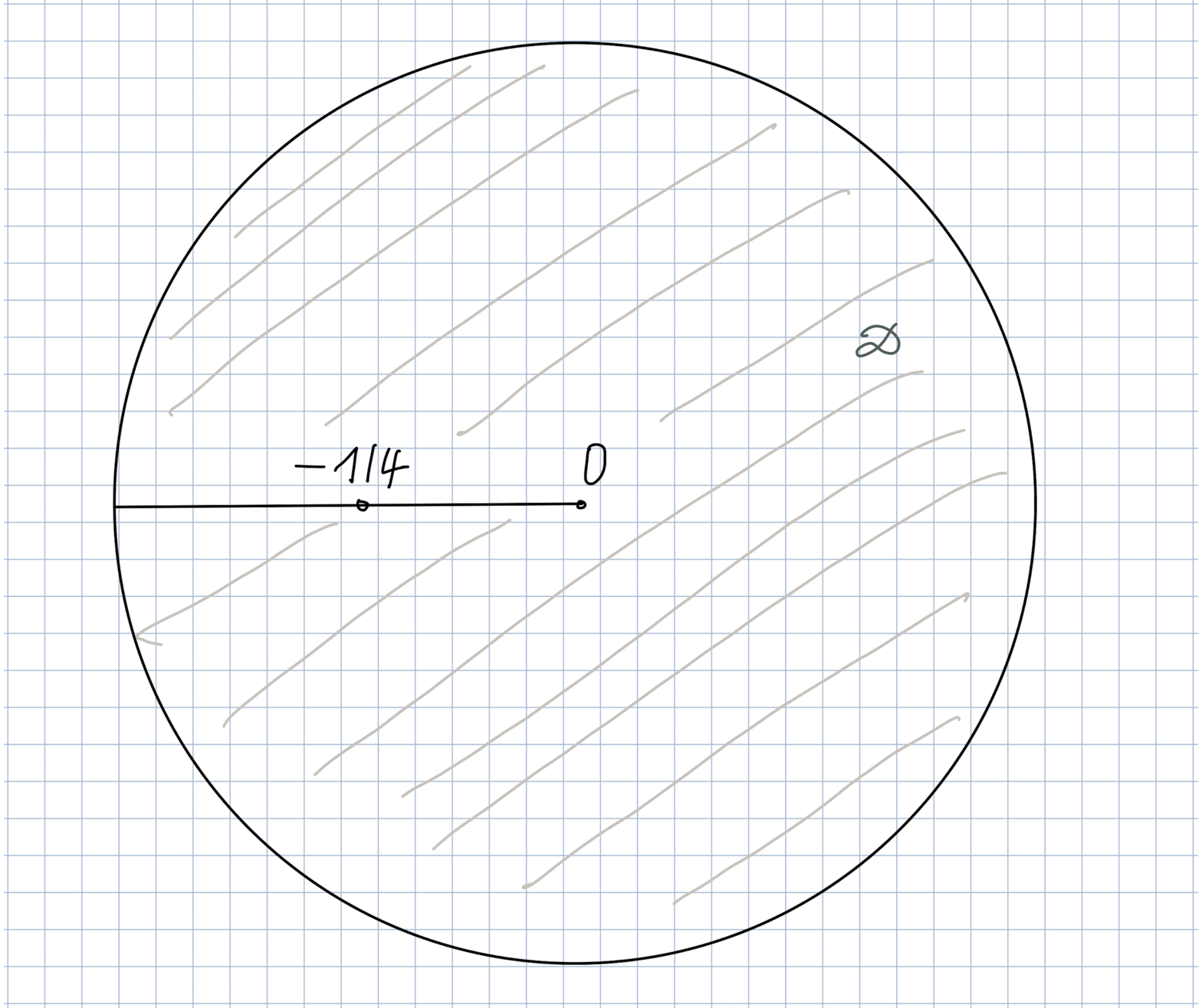}
\end{center}
\caption{ The complex domain $\mathcal D = \{ \C \setminus (-\infty, 0] \}$}
\label{figure2a}
\end{figure}
and hence the complete elliptic integral $J(h)$ allows an analytic continuation on this domain too, see Fig. \ref{figure2a}.
As $J_0$ is a complete elliptic integral of first kind, it is a period of $E_h$, so it can not vanish for $ h\in \mathcal D$.  Therefore we consider the function
\begin{align}
\label{fh}
F(h) = \frac{J(h)}{J_0(h)} = \tilde p(h) \frac{J_2(h)}{J_0(h)} +  \tilde q(h) \;
, h\in \mathcal D
\end{align}
and we shall equivalently count its zeros in  $ \mathcal D$. 
The cycle $\gamma(h)$ has a limit when $h$ tends to the half-open interval $(-\infty, 0]$ and we denote
the resulting cycle by $\gamma (h^+)$ or $\gamma(h^-)$, depending on 
whether $Im(h) >0$, or  $Im(h)<0$, throughout the limit process. 

For real $h$ the curve $E_h$ is real and allows an anti-holomorphic involution (complex conjugation) which induces an involution on $ H_1(E_h, \Z)$. It is easy to see that
$$
\overline{\gamma(h^+)} = \gamma(h^-) , h\leq 0
$$
and that $\gamma(h^+)$ is a continuation of $\gamma(h^-)$ along a path contained in $\mathcal D$. Therefore $\gamma(h^+)$ is computed from $\gamma(h^-)$ by the Picard-Lefschetz formula, involving the cycles vanishing at the critical values $h=0$ and $h=-\frac14$ of $H$. More precisely, denote for $-\frac 14 < h< 0$ by $\gamma_1(h), \gamma_2(h)$ the two real ovals 
of the curve $E_h$, and denote by $\delta(h)$ the cycle vanishing at $(0,0)$ when $h\to 0$. By the Picard-Lefschetz formula we have
\begin{align*}
\gamma(h^+)-\gamma(h^-) &= 2\delta(h) \\
\gamma(h^+) &= \gamma_1(h) + \gamma_2(h) + \delta(h), \gamma(h^-) = \gamma_1(h) + \gamma_2(h) - \delta(h) 
\end{align*}
and in particular $$
\overline{\delta(h) } = - \delta(h) , \overline{\gamma_1(h) } = \gamma_1(h),  \overline{\gamma_2(h) }= \gamma_2(h).
$$
Clearly in $H_1(E_h,\Z)$ holds $\gamma_1(h)=\gamma_2(h)$, and the cycles $\gamma(h^+) $,  $\gamma(h^-) $, $-\frac14 < h < 0$, generate $H_1(E_h,\Z)$.
In fact  $\gamma(h^+) $,  $\gamma(h^-) $ generate  $H_1(E_h,\Z)$ also for 
$h \in (-\infty, -\frac14) $ (again by the Picard-Lefschetz formula).

For two cycles $\alpha, \beta \in H_1(E_h,\Z)$ consider the Wronskian
$$
W(\alpha,\beta) =  \det \left(\begin{array}{cc}\int_{\alpha} \frac{dx}{y} & \int_{\alpha} \frac{x^2dx}{y} \\\int_{\beta} \frac{dx}{y} & \int_{\beta} \frac{x^2dx}{y} \end{array}\right) 
$$
and denote
\begin{align}
\label{w1}
W_1 & = W(\gamma(h^+),\gamma(h^-)), -\frac14 < h < 0 \\
\label{w2}
W_2 &= W(\gamma(h^+),\gamma(h^-)), h < -\frac14 .
\end{align}
Let $\alpha(h), \beta(h) $ be continuous families of cycles which generate a basis of $H_1(E_h,\Z)$.
The Picard-Lefschetz formula and the moderate growth of the corresponding  integrals imply that $W(\alpha(h),\beta(h))$ 
is a rational function  in $h$. 
which has neither zeros, nor poles on the complex plane $\C$, except eventually at $h=0, -\frac 14$. A local analysis near $0, -\frac14, \infty$ then shows that $W(h)$ has no poles at all, and hence is a non-zero imaginary constant. In particular, the Wronskians $W_1, W_2$ defined in (\ref{w1}), (\ref{w2}) are non-zero imaginary constants.
A crucial point in our proof will be the following observation
\begin{lemma}
\label{lemma4}
$W_1/W_2  $ is a real strictly positive constant.
\end{lemma}
\proof
On the interval $(-\frac14,0)$ we have
\begin{align*}
W_1 = -2W( \gamma_1(h) + \gamma_2(h), \delta(h) ) .
\end{align*}
To compare $W_1$ to $W_2$ we consider the analytic continuation of 
$W_2= W(\gamma(h^+),\gamma(h^-))$ near a point $h\in(-\infty,-\frac 14)$, along a path in the upper half-plane $Im(h)>0$ to a function 
near $h$ in the interval $(-\frac14,0)$. The initial cycle $\gamma(h^+)$ is deformed to
$$
\gamma(h^+) = \gamma_1(h) + \gamma_2(h) + \delta(h) ,  -\frac14 < h < 0
$$
and the initial cycle $\gamma(h^-)$ is deformed, according to the Picard-Lefschetz formula, to the cycle
$$
\gamma(h^-)  \pm ( \gamma_1(h) + \gamma_2(h)),  -\frac14 < h < 0 .
$$
Therefore
\begin{align*}
W_2 &= W_1 \pm W(\gamma_1(h) + \gamma_2(h) + \delta(h), \gamma_1(h) + \gamma_2(h))\\
&=  W_1\pm W(\delta(h), \gamma_1(h) + \gamma_2(h)) = W_1\pm \frac12 W_1 
\end{align*}
which completes the proof of the Lemma. \endproof
We are ready to compute the zeros of $F$, see (\ref{fh}), in the  complex domain $\mathcal D$. 
We apply the argument principle to the domain bordered  by the cut $(-\infty,0]$ and a  circle of a sufficiently big radius.

\begin{itemize}
\item 
The function $F$ has a real limit at $h=0$.
\item
Along the cut $(-\infty, 0) $ the imaginary part of $F$ equals
$$
Im (F(h)) = \tilde p(h)\frac{W(\gamma(h^+),\gamma(h^-))}{|J_0(h)|^2}
$$
where the Wronskian $W(\gamma(h^+),\gamma(h^-))$ is a non zero constant on the interval $(-\infty, - 1/4) $, a 
non zero constant on $(-1/4,0)$, and its limit at $h=-\frac14$  is $0$. Therefore along  $(-\infty,0)$  the imaginary part of $F$
vanishes at most three times, at $h=-\frac14$ and the eventual zeros of  $\tilde p(h) $. Therefore the argument of $F(h)$ increases by at most $6\pi$. By Lemma \ref{lemma4}, however, 
crossing the point $h=-\frac14$ does not contribute to the increase of the argument of $F$. Thus, the argument of $F(h)$ along  $(-\infty,0)$ increases by at most $4\pi$. 
\item
Along a circle with sufficiently big radius, the functions $F(h)$ behaves like
$const. \times h^\alpha$ where $0 \leq \alpha \leq  2+1/2 $. The increase of the argument of $F$ along such a  circle is close to $5\pi$ or less. 
\end{itemize}
Summing up the above information, we conclude that the increase of the argument of $F$ along the border of the domain is at most $10\pi$. Thus $F(h)$ can have at most five zeros, counted with multiplicity, in the complex domain $\mathcal D$. This implies that 
$J(h)=I'(h)$ has at most five zeros in $(0,\infty)$ and hence $I(h)$ has at most six zeros  in $(0,\infty)$.
This completes the study of the exterior period annulus.

The study of the interior period annulus of the eight-loop is analogous, in the complex domain $\mathcal D = \C\setminus [0,\infty)$, with the same conclusion : $J$ has at most five zeros, and $I$ has at most six zeros on $[-\frac14,0)$. However, $I(-\frac14)=0$ so finally $I$ has at most five zeros in $(-\frac14,0)$. 

This result can be further improved, for the interval  $(-\frac14,0]$, where $h=0$ corresponds to a saddle loop connection of the non-perturbed system. For this purpose we define 
the cyclicity of $h=0$ as the maximal $k$ such that 
$$
J(h)= c_0+c_1 h 
$$

Similar result with the same proof holds for the truncated pendulum and the global center. 
Theorem \ref{th6} is proved.

\section{The Nakayama lemma and the Bautin ideal}

Let $K\{\lambda\} = K\{\lambda_1,\dots,\lambda_n\}$ be the ring of  germs of analytic functions at the origin in 
$\lambda$ where $K= \R, \C$,  $\mathfrak m= (\lambda_1,\lambda_2,\dots,\lambda_n)$ be the maximal ideal of $K\{\lambda\}$.
Let $$b_i=b_i^0+b_i^1 \in K\{\lambda\}  , \;\; i=1,2,\dots, k$$ 
be germs of analytic functions vanishing at the origin and
consider the ideals
$$
\mathcal B = (b_1,b_2, \dots,b_k), \; \mathcal B^0 = (b_1^0,b_2^0, \dots,b_k^0), \;  \mathcal B^1 = (b_1^1,b_2^1, \dots,b_k^1) .
$$

The next result is a version of
Nakayama lemma from Algebraic Geometry, see \cite[Lemma 7.4]{bry10}.

\begin{lemma}
\label{naka}
If $\mathcal B^1 \subset \mathfrak m\, \mathcal B^0$ then $\mathcal B = \mathcal B^0$.
\end{lemma}
\begin{example}
Let
\begin{align*}
b_1= \lambda_1^2 +  \lambda_1^2\lambda_2^2+  \lambda_1\lambda_2^3 + \lambda_2^4 , \; b_2 = \lambda_2^3 + \lambda_1^4 + \lambda_1^3\lambda_2
\end{align*}
The polynomials $b_1,b_2$ belong to the   ideal generated by $\lambda_1^2, \lambda_2^3$,
$$
(b_1,b_2) \subset (\lambda_1^2, \lambda_2^3) .
$$
By the above Lemma in the local ring $ K\{\lambda\} $ holds 
$$
(b_1,b_2) = (\lambda_1^2, \lambda_2^3),
$$
that is to say $\lambda_1^2, \lambda_2^3$ are linear combinations of $b_1,b_2$ whose  coefficients are suitable
analytic functions in  the local ring $ K\{\lambda\} $.
\end{example}
The above turns out to be a useful tool, for determining  the generators of the localized Bautin ideal. Indeed, if for some reasons the Bautin ideal is generated by $b_1,b_2$ as above, then it is also locally generated by $\lambda_1^2, \lambda_2^3$ which is much simpler.\\

\proof[Proof of Lemma \ref{naka}]
Assume that  $\mathcal B^1 \subset \mathfrak m\, \mathcal B^0$, that is to say
\begin{align}
\label{bn}
 b_i^1 =\sum_{j=1}^k a_{ij} b_j^0 , \; i=1,2,\dots,k 
\end{align}
where $a_j= a_j(\lambda) \in  \mathfrak m $ are germs of analytic functions vanishing at the origin.
The relations  (\ref{bn}) can be written in a matrix form
$$
\left(\begin{array}{c} \\b_1 \\b_2 \\ \vdots \\ b_k\end{array}\right)
= \left(\begin{array}{cccc}1+a_{11} & a_{12} & \dots & a_{1k} \\
 a_{21} &1+a_{22} & \dots & a_{2k} \\ \vdots&  \vdots &  \vdots &  \vdots \\a_{k1} & a_{k_2} & \dots & 1+a_{kk}\end{array}\right)
 \left(\begin{array}{c} \\b_1^0 \\b_2^0 \\ \vdots \\ b_k^0 \end{array}\right)
$$
ot simply
\begin{align}
\label{bn1}
b= (I+A)b^0
\end{align}
where $b, b^0$ are vector-columns with entries $b_i$, $b_i^0$ respectively, $I$ is the identity matrix, and
$A= (a_{ij})$. Obviously, the ideal $(b_1^0,b_2^0,\dots,b_k^0)$ contains the ideal
$(b_1,b_2, \dots,b_k)$.
We use the Cramer's rule to write
$$
b^0 = (I+A)^{-1} b
$$
in the form
\begin{align}
\label{bn2}
b_i^0 = b_i + \sum_{j=1}^k b_j \tilde a_{ij} , \; i=1,2,\dots,k 
\end{align}
where $\tilde a_{ij}$ are germs of analytic functions vanishing at the origin in the variables $(b_1^0,b_2^0,\dots,b_k^0)$. This shows that the ideal $(b_1,b_2, \dots,b_k)$ contains the ideal
 $(b_1^0,b_2^0,\dots,b_k^0)$.
\endproof

The Nakayama lemma turns out to be an important tool in the study of the localized Bautin ideal.
To illustrate  this  we shall describe in detail the localized Bautin ideal of the perturbed vector field $X_\lambda$, see (\ref{eq2}), in the
 eight loop case (iii).

We consider, for small parameters $\lambda$, the first return map and denote by $\Delta_\lambda(h)$ the corresponding displacement map. 
As usual $h$ is a local variable on a cross-section to the orbits and containing the equilibrium point. The displacement map is analytic both in $\lambda$ and $h$, and we may expand
$$
d(h,\lambda)= \sum_{i\geq 1} a_i  h^i
$$
where the analytic function $a_i=a_i(\lambda)$ vanish at the origin and generate an ideal in the local ring $K\{\lambda\}$. A basic fact about this localized Bautin ideal is that it can be studied via the corresponding arc space, that is to say one-parameter space of deformations of the vector field, which on its turn allows to use the powerful machinery of Melnikov functions, e.g. \cite{gavr20}.
Namely, it follows from Theorem \ref{th4} that the linear span of all Melnikov function $M_k$ is a six-dimensional vector space. The localized Bautin ideal has therefore at least six generators. We denote these generators $b_i$ and write
$$
d(h,\lambda)= \sum _{i=1}^N b_i(\lambda) \varphi_i(h,\lambda) , N\geq 6
$$
where the functions $\varphi_i(h,0)$ 	are linearly independent. We may assume moreover that, $i=1, 2, \dots, N$ are linearly independent Melnikov functions, enumerated in section \ref{section12}. 
This observation allows to reconstruct explicitly
the generators $b_i$, $i\leq 6$, and subsequently show that $N=6$.

Namely, consider instead of a general 6-parameter deformation, just an arc (one-parameter deformation)
$$
\varepsilon \to \lambda(\varepsilon) 
$$ 
in the parameter space. The expansion for $d(h,\lambda)$ takes the form
$$
d(h,\lambda(\varepsilon))(h) = \varepsilon^k M_k(h) + \dots
$$
where $M_k$ are linear combinations of $\varphi_1(h,0),\dots, \varphi_6(h,0)$.
More explicitly, take a deformation  linear in $\varepsilon$
$$\varepsilon \to \lambda =   (\varepsilon \lambda_{11}, \varepsilon \lambda_{21}, \dots, \varepsilon \lambda_{61}) .$$
The linear (or first) Melnikov functions provide the order one approximation in $\lambda$  of the displacement map.
Indeed, it follows from Theorem \ref{th4}, that $M_1=0$ if and only if
\begin{align}
\label{b1}
\lambda_{11} = \lambda_{21} - 3\lambda_{31} = \lambda_{41} + 3 \lambda_{31} = \lambda_{51} = \lambda_{61}=0
\end{align}
which implies that the first order approximations of the six generators $b_i$ are given by the above five linear functions. Therefore, we may assume that
\begin{align*}
b_1&=\lambda_{1}+ \dots,  b_2 = \lambda_{2}- 3\lambda_{3} + \dots,  b_3= \dots \\
b_4&= \lambda_{4} + 3 \lambda_{3} + \dots,   b_5= \lambda_{5} + \dots,  b_6= \lambda_{6} + \dots
\end{align*}
where the dots replace some higher order terms. For a further use, let us denote the linear terms 
\begin{align*}
b_1^0&=\lambda_{1},  b_2^0 = \lambda_{2}- 3\lambda_{3} \\
b_4^0&= \lambda_{4} + 3 \lambda_{3},   b_5^0= \lambda_{5},  b_6^0 = \lambda_{6} .
\end{align*}
To discover the second order terms of $b_i$, $i\leq 6$, we assume that 
(\ref{b1}) holds true and compute the second Melnikov function $M_2$, which according to Theorem \ref{th4} vanishes identically. This means that all second order terms of the generators $b_i$ are in the ideal 
$$
\tilde{\mathcal B } = (b_1^0, b_2^0, b_4^0, b_5^0, b_6^0 ) .
$$
and as they are quadratic, they also belong to 
$
 \mathfrak m\, \tilde{\mathcal B } .
$
We proceed therefore to the computation of the third order Melnikov function $M_3$. 
It is seen from Theorem \ref{th4}, that the function $M_3$ does not belong to the five-dimensional vector space formed by the first Melnikov functions, and moreover it is multiplied by $\lambda_{31}^3$. Therefore we may put
$$
b_3= \lambda_3^3 + \dots
$$
where the dots replace terms which either belong to $ \mathfrak m\, \tilde{\mathcal B } $ or they are
of order at least four, that is to say of the form $\lambda_3^k, k\geq 4$.  To resume, if we put
$$
b_3^0= \lambda_3^3
$$
then we have
$$
b_i= b_i^0 + b_i^1, \;\; i=1, \dots, 6 
$$
where $b_i^1$ belong to the ideal  $\mathfrak m\, \mathcal B _0 $ where 
$$
\mathcal B _0  = (b_1^0, b_2^0, b_3^0, b_4^0, b_5^0, b_6^0 ) 
$$
and 
$$
\mathcal B   = (b_1 , b_2 , b_3 , b_4 , b_5 , b_6 ) .
$$
According to Lemma \ref{naka} we have $\mathcal B = \mathcal B_0 $. As for $b_j$, $j>6$ we note that these terms belong to the radical of $\mathcal B$. Indeed, if some $b_j$ is not in the radical of $\mathcal B$, then there would be an arc for which
$\varphi_j(h,0)$ is a new Melnikov function, which is not the case. But if  $b_j$ is  in the radical of $\mathcal B$, then it is in the ideal generated by $\mathcal B, \lambda_3^k$ where $k< 3$. But in this case, once again we can find an arc, for which $\varphi_j(h,0)$ is a Melnikov function which can not be the case. In conclusion, $b_j$, $j>6$ are well in $\mathcal B$ and the Theorem \ref{th5} (iii) is proved. Finally, the six-parameter deformations of the global center and the truncated pendulum are studied in the same way. 
$\Box$

We note that in the eight loop case we have three period annuli but their Bautin ideals are the same.

\section{Limit Cycles}
In this section, following \cite{ilie98,gavr08,gavr20}, we determine the exact upper bounds for the number of limit cycles, bifurcating from a given period annulus of the 6-parameter perturbed equation (\ref{eq1}). Indeed, according to the preceding section, the displacement map allows an expansion
\begin{align*}
d(h,\lambda)= \sum _{i=1}^6 b_i(\lambda) (\varphi_i(h,0)  + O(\lambda)).
\end{align*}
where $\varphi_i(h,0)$ form a basis of the 6-dimensional vector space of Abelian integrals
\begin{equation}
\label{i02}
V= Span\{ h^iI_0(h), h^jI_2(h): \; i,j = 0,1, 2 \}
\end{equation}
described in  Theorems \ref{th1}-\ref{th4}, and $b_i(\lambda)$ are the generators of the Bautin ideal $\mathcal B$
described in  Theorem \ref{th5}. As $\mathcal B$ is not radical we consider its blow up \cite[section 2]{gavr20}
and observe that the exceptional divisor over $\lambda=0$ is just the projective space $\P^5$. Thus not only every Melnikov function belongs to the space V in (\ref{i02}) 
but every function in $V$ is realized as a Melnikov function of an appropriate one-parameter deformation.

As the ideal sheaf defined by $\mathcal B$ along the exceptional divisor is locally principal,
we may use arcs (one-parameter deformations) to study the cyclicity of the period annuli, see  \cite[Theorem 1]{gavr08}. 
Therefore   Theorem \ref{th6} implies   Theorem \ref{main}.


\section{Appendix}

Below we list the relative decomposition formulas we used in the global center case $H=\frac12y^2+\frac12x^2+\frac14x^4$.  
 
$$
\begin{array}{ll}

(i1)    &  y^3dx=(\frac{12}{7}H-\frac37x^2)ydx-\frac37xydH+d(\frac17xy^3)\\[2mm]
(i2)    & x^4ydx=(\frac47H-\frac87x^2)ydx+\frac67xydH-d(\frac27xy^3)\\[2mm]
(i3)    & y^5dx= [\frac{240}{77}H^2-\frac{320}{231}Hx^2-\frac{40}{231}x^2+\frac{20}{231}H]ydx\\[2mm]
          & +[\frac{10}{77}xy+\frac{5}{33}x^3y-\frac{60}{77}Hxy-\frac57xy^3]dH\\[2mm]
          & +d[\frac{20}{77}Hxy^3+\frac{1}{11}xy^5-\frac{10}{231}xy^3-\frac{5}{99}x^3y^3]\\[2mm]
 (i4)   & x^6ydx=(\frac43Hx^2+\frac{32}{21}x^2-\frac{16}{21}H)ydx+(\frac23x^3y-\frac87xy)dH
             -d(\frac29x^3y^3-\frac{8}{21}xy^3)\\[2mm]

 (i5)    & x^2y^3dx=(\frac43Hx^2+\frac{8}{21}x^2-\frac{4}{21}H)ydx-(\frac27xy+\frac13x^3y)dH
         +d(\frac19x^3y^3+\frac{2}{21}xy^3)\\[2mm]
 (i6)   & x^2y^5dx =(\frac{80}{39}H^2x^2+\frac{3620}{3003}Hx^2-\frac{1600}{3003}H^2-\frac{80}{1001}H
              +\frac{160}{1001}x^2)ydx\\[2mm]
          & -(\frac{380}{1001}Hxy+\frac{20}{39}Hx^3y+\frac{130}{273}xy^3+\frac59x^3y^3+\frac{20}{143}x^3y
             +\frac{120}{1001}xy)dH\\[2mm]
          & +d(\frac{1}{13}x^3y^5+\frac{10}{143}xy^5+\frac{20}{91}Hx^3y^3+\frac{20}{429}x^3y^3+\frac{40}{1001}xy^3
              +\frac{380}{3003}Hxy^3)\\[2mm]

(i7)  & xy^4dx=-\frac23(4Hx^2+x^2+y^2+x^2y^2)dH\\[2mm]
      & +d[2x^2H^2+\frac23H^2-Hx^4+(\frac16-\frac13H)x^6+\frac18x^8+\frac{1}{40}x^{10}]\\[2mm]
(i8) & x^ky^2dx=-\frac{2}{k+1}x^{k+1}dH+d[\frac{2H}{k+1}x^{k+1}-\frac{1}{k+3}x^{k+3}-\frac{1}{2(k+5)}x^{k+5}]\\[3mm]
(i9) & x^ky^4dx=-(\frac{8H}{k+1}x^{k+1}-\frac{4}{k+3}x^{k+3}-\frac{2}{k+5}x^{k+5})dH\\[2mm]
     & \hspace{13mm} +d[\frac{4H^2}{k+1}-\frac{4H}{k+3}x^2+\frac{1-2H}{k+5}x^4+\frac{1}{k+7}x^6
        +\frac{1}{4(k+9)}x^8]x^{k+1}

\end{array}
$$
\newpage
The following are the relative decomposition formulas in the truncated pendulum case $H=\frac12y^2+\frac12x^2-\frac14x^4$.

$$\begin{array}{ll}
(ii1)    &  y^3dx=(\frac{12}{7}H-\frac37x^2)ydx-\frac37xydH+d(\frac17xy^3)\\[2mm]
(ii2)    & x^4ydx=(\frac87x^2-\frac47H)ydx-\frac67xydH+d(\frac27xy^3)\\[2mm]
(ii3)    & y^5dx= [\frac{240}{77}H^2-\frac{320}{231}Hx^2+\frac{40}{231}x^2-\frac{20}{231}H]ydx\\[2mm]
          & -[\frac{10}{77}xy-\frac{5}{33}x^3y+\frac{60}{77}Hxy+\frac57xy^3]dH\\[2mm]
          & +d[\frac{20}{77}Hxy^3+\frac{1}{11}xy^5+\frac{10}{231}xy^3-\frac{5}{99}x^3y^3]\\[2mm]
 (ii4)   & x^6ydx=-(\frac43Hx^2-\frac{32}{21}x^2+\frac{16}{21}H)ydx-(\frac23x^3y+\frac87xy)dH
             +d(\frac29x^3y^3+\frac{8}{21}xy^3)\\[2mm]

 (ii5)    & x^2y^3dx=(\frac43Hx^2-\frac{8}{21}x^2+\frac{4}{21}H)ydx+(\frac27xy-\frac13x^3y)dH
         +d(\frac19x^3y^3-\frac{2}{21}xy^3)\\[2mm]
 (ii6)   & x^2y^5dx =(\frac{80}{39}H^2x^2-\frac{3620}{3003}Hx^2+\frac{1600}{3003}H^2-\frac{80}{1001}H
              +\frac{160}{1001}x^2)ydx\\[2mm]
          & +(\frac{380}{1001}Hxy-\frac{20}{39}Hx^3y+\frac{130}{273}xy^3-\frac59x^3y^3+\frac{20}{143}x^3y
             -\frac{120}{1001}xy)dH\\[2mm]
          & +d(\frac{1}{13}x^3y^5-\frac{10}{143}xy^5+\frac{20}{91}Hx^3y^3-\frac{20}{429}x^3y^3+\frac{40}{1001}xy^3
              -\frac{380}{3003}Hxy^3)\\[2mm]

(ii7)  & xy^4dx=-\frac23(4Hx^2-x^2-y^2+x^2y^2)dH\\[2mm]
      & +d[2x^2H^2-\frac23H^2-Hx^4+(\frac16+\frac13H)x^6-\frac18x^8+\frac{1}{40}x^{10}]\\[2mm]
\end{array}$$

\noindent
The relative decomposition formulas in the eight loop case $H=\frac12y^2-\frac12x^2+\frac14x^4$, of them
$(iii1), (iii2), (iii5)$ are taken from  \cite{ILY}.  

$$\begin{array}{ll}
(iii1)    &  y^3dx=(\frac{12}{7}H+\frac37x^2)ydx-\frac37xydH+d(\frac17xy^3)\\[2mm]
(iii2)    & x^4ydx=(\frac47H+\frac87x^2)ydx+\frac67xydH-d(\frac27xy^3)\\[2mm]
(iii3)    & y^5dx= [\frac{240}{77}H^2+\frac{320}{231}Hx^2+\frac{40}{231}x^2+\frac{20}{231}H]ydx\\[2mm]
          & +[\frac{10}{77}xy-\frac{5}{33}x^3y-\frac{60}{77}Hxy-\frac57xy^3]dH\\[2mm]
          & +d[\frac{20}{77}Hxy^3+\frac{1}{11}xy^5-\frac{10}{231}xy^3+\frac{5}{99}x^3y^3]\\[2mm]
(iii4)    & x^6ydx=(\frac43Hx^2+\frac{32}{21}x^2+\frac{16}{21}H)ydx+(\frac23x^3y+\frac87xy)dH
             -d(\frac29x^3y^3+\frac{8}{21}xy^3)\\[2mm]

(iii5)    & x^2y^3dx=(\frac43Hx^2+\frac{8}{21}x^2+\frac{4}{21}H)ydx+(\frac27xy-\frac13x^3y)dH
         +d(\frac19x^3y^3-\frac{2}{21}xy^3)\\[2mm]
(iii6)    & x^2y^5dx =(\frac{80}{39}H^2x^2+\frac{3620}{3003}Hx^2+\frac{1600}{3003}H^2+\frac{80}{1001}H
              +\frac{160}{1001}x^2)ydx\\[2mm]
          & +(\frac{380}{1001}Hxy-\frac{20}{39}Hx^3y+\frac{130}{273}xy^3-\frac59x^3y^3-\frac{20}{143}x^3y
             +\frac{120}{1001}xy)dH\\[2mm]
          & +d(\frac{1}{13}x^3y^5-\frac{10}{143}xy^5+\frac{20}{91}Hx^3y^3+\frac{20}{429}x^3y^3-\frac{40}{1001}xy^3
              -\frac{380}{3003}Hxy^3)\\[2mm]

(iii7)  & xy^4dx=-\frac23(4Hx^2+x^2-y^2+x^2y^2)dH\\[2mm]
      & +d[2x^2H^2-\frac23H^2+Hx^4+(\frac16-\frac13H)x^6-\frac18x^8+\frac{1}{40}x^{10}]\\[2mm]
\end{array}$$

\end{document}